\documentclass[12pt]{article}
\usepackage{amssymb,amsmath,mathrsfs}
\usepackage{hyperref}
\usepackage{graphicx}
\usepackage{color}
\usepackage[OT2,OT1]{fontenc}
\usepackage{upquote}

\begin{document}

\title{\LARGE\bf On the equivalence of Fourier expansion and Poisson summation formula for the series approximation of the exponential function}

\author{
\normalsize\bf S. M. Abrarov\footnote{\scriptsize{Dept. Earth and Space Science and Engineering, York University, Toronto, Canada, M3J 1P3.}}\,, B. M. Quine$^{*}$\footnote{\scriptsize{Dept. Physics and Astronomy, York University, Toronto, Canada, M3J 1P3.}} and R. K. Jagpal$^\dagger$}

\date{April 25, 2019}
\maketitle

\begin{abstract}
In this short note we show the equivalence of Fourier expansion and Poisson summation approaches for the series approximation of the exponential function $\exp \left( { - {t^2}/4} \right)$. The application of the Poisson summation formula is shown to reduce to that of the Fourier expansion method.
\end{abstract}

\section {Methodology}

Let us apply the Poisson summation formula (see for example \cite{Beerends2003}) to the exponential function $\exp \left( { - {t^2}/4} \right) \equiv \exp \left[ { - {{\left( {t/2} \right)}^2}} \right]$, i.e.:
$$
\sum\limits_{n =  - \infty }^\infty  {\exp \left[ {{{-\left( {  \frac{t}{2} + n{\tau _m}} \right)}^2}} \right]}  = \frac{{\sqrt \pi  }}{{{\tau _m}}}\sum\limits_{n =  - \infty }^\infty  {\exp \left( { - \frac{{{\pi ^2}{n^2}}}{{\tau _m^2}}} \right)} \cos \left( {\frac{{\pi n}}{{{\tau _m}}}t} \right)
$$
or
\small
\begin{equation}\label{eq_1}
\begin{aligned}
\exp \left( {- \frac{{{t^2}}}{4}} \right) &= \frac{{\sqrt \pi  }}{{{\tau _m}}}\sum\limits_{n =  - \infty }^\infty  {\exp \left( { - \frac{{{\pi ^2}{n^2}}}{{\tau _m^2}}} \right)} \cos \left( {\frac{{\pi n}}{{{\tau _m}}}t} \right) \\ 
&  - \left\{ {\sum\limits_{n =  - \infty }^{ - 1} {\exp \left[ {{{-\left( {  \frac{t}{2} + n{\tau _m}} \right)}^2}} \right]}  + \sum\limits_{n = 1}^\infty  {\exp \left[ {{{-\left( {  \frac{t}{2} + n{\tau _m}} \right)}^2}} \right]} } \right\}, \\
\end{aligned}
\end{equation}
\normalsize
where $2{\tau _m}$ is the period. Assuming that the half-period is large enough, say ${\tau _m} \geqslant 12$, within the range $t \in \left[ { - {\tau _m},{\tau _m}} \right]$ it follows that
\scriptsize
\begin{equation*}
\exp {\left( { - \frac{t^2}{4}} \right)} >  > \exp \left[ {{{-\left( {  \frac{t}{2} \pm {\tau _m}} \right)}^2}} \right] >  > \exp \left[ {{{-\left( {  \frac{t}{2} \pm 2{\tau _m}} \right)}^2}} \right] >  > \exp \left[ {{{-\left( {  \frac{t}{2} \pm 3{\tau _m}} \right)}^2}} \right] ...\, .
\end{equation*}
\normalsize
Consequently, we can ignore the terms in curly brackets in equation \eqref{eq_1} and approximate exponential function as
\begin{equation}\label{eq_2}
\exp \left( { - \frac{{{t^2}}}{4}} \right) \approx \frac{{\sqrt \pi  }}{{{\tau _m}}}\sum\limits_{n =  - \infty }^\infty  {\exp \left( { - \frac{{{\pi ^2}{n^2}}}{{\tau _m^2}}} \right)} \cos \left( {\frac{{\pi n}}{{{\tau _m}}}t} \right).
\end{equation}
Taking into account that
$$
\exp \left( { - \frac{{{\pi ^2}{{\left( { - n} \right)}^2}}}{{\tau _m^2}}} \right)\cos \left( {\frac{{\pi \left( { - n} \right)}}{{{\tau _m}}}t} \right) \equiv \exp \left( { - \frac{{{\pi ^2}{n^2}}}{{\tau _m^2}}} \right)\cos \left( {\frac{{\pi n}}{{{\tau _m}}}t} \right),
$$
the approximation \eqref{eq_2} can be rearranged in form
$$
\exp \left( { - \frac{{{t^2}}}{4}} \right) \approx \frac{{\sqrt \pi  }}{{{\tau _m}}} + \frac{{2\sqrt \pi  }}{{{\tau _m}}}\sum\limits_{n = 1}^\infty  {\exp \left( { - \frac{{{\pi ^2}{n^2}}}{{\tau _m^2}}} \right)} \cos \left( {\frac{{\pi n}}{{{\tau _m}}}t} \right)
$$
or
\begin{equation*}
\exp \left( { - \frac{{{t^2}}}{4}} \right) \approx  - \frac{{\sqrt \pi  }}{{{\tau _m}}} + \frac{{2\sqrt \pi  }}{{{\tau _m}}}\sum\limits_{n = 0}^\infty  {\exp \left( { - \frac{{{\pi ^2}{n^2}}}{{\tau _m^2}}} \right)} \cos \left( {\frac{{\pi n}}{{{\tau _m}}}t} \right).
\end{equation*}
Lastly, defining
$$
{a_n} = \frac{{2\sqrt \pi  }}{{{\tau _m}}}\exp \left( { - \frac{{{\pi ^2}{n^2}}}{{\tau _m^2}}} \right)
$$
we obtain
\begin{equation}\label{eq_3}
\exp \left( { - \frac{t^2}{4}} \right) \approx  - \frac{{{a_0}}}{2} + \sum\limits_{n = 0}^\infty  {{a_n}} \cos \left( {\frac{{\pi n}}{{{\tau _m}}}t} \right),\qquad t \in \left[ { - {\tau _m},{\tau _m}} \right].
\end{equation}

We derived equation \eqref{eq_3} in a different way by Fourier series expansion and used it for numerical integration of the Voigt function by truncating its upper integration limit such that \cite{Abrarov2010}
\[
\begin{aligned}
K\left(x,y\right)&=\frac{1}{\sqrt{\pi}}\int_0^\infty{\exp\left(-t^2/4\right)\exp\left(-yt\right)\cos\left(xt\right)}\\
&\approx \frac{1}{\sqrt{\pi}}\int_0^{\tau_m}{\exp\left(-t^2/4\right)\exp\left(-yt\right)\cos\left(xt\right)},
\end{aligned}
\]
where $x$ and $y$ are input parameters.

It should be noted that application of the Poisson summation formula to the functions of kind $\exp\left(-t^2\right)$
and $\exp\left(-t^2\right)/(t-\alpha)$ such that $\alpha\in\mathbb{C}\backslash\hspace{-3pt}\left\{0\right\},$
is a very efficient method in numerical integrations that was used, for example, in approximations of the Dawson\text{'}s integral \cite{Dawson1897}
$$
\text{daw}\left(z\right)= \exp\left(-z^2\right)\int_0^z\exp\left(t^2\right) dt
$$
and the error function \cite{Salzer1951, Molin2011}
$$
\text{erf}\left(z\right) = \frac{2}{\sqrt{\pi}}\int_0^z\exp\left(-t^2\right) dt,
$$
where $z = x+iy$ is the complex argument.

\section{Conclusion}

We show the equivalence of series approximation derived by Fourier expansion and by Poisson summation formula for the exponential function $\exp \left( { - {t^2}/4} \right)$ that we applied for numerical integration of the Voigt function \cite{Abrarov2010}.

\bigskip

\end{document}